\documentclass[%
 amsmath,amssymb,
prfluids
]{revtex4-2}

\usepackage{graphicx}
\usepackage{dcolumn}
\usepackage{bm}

\usepackage{multirow}
\usepackage{array}

\usepackage{algorithm}
\usepackage{algpseudocode}

\renewcommand{\b}[1]{\boldsymbol{#1}}
\newcommand{\T}[0]{^\mathrm{T}}
\newcommand{\tr}[0]{\mathrm{tr}}
\newcommand{\att}[0]{\mathrm{att}}

\begin{document}

\title{Multi-objective SINDy for parameterized model discovery from single transient trajectory data}

\author{Javier A. Lemus}
\author{Benjamin Herrmann}%
 \email{benjaminh@uchile.cl}
\affiliation{Department of Mechanical Engineering, University of Chile, Beauchef 851, Santiago, Chile}


\begin{abstract}
The sparse identification of nonlinear dynamics (SINDy)~\citep{brunton2016pnas} has been established as an effective technique to produce interpretable models of dynamical systems from time-resolved state data via sparse regression. %
However, to model parameterized systems, SINDy requires data from transient trajectories for various parameter values over the range of interest, which are typically difficult to acquire experimentally.
In this work, we extend SINDy to be able to leverage data on fixed points and/or limit cycles to reduce the number of transient trajectories needed for successful system identification.
To achieve this, we incorporate the data on these attractors at various parameter values as constraints in the optimization problem.
First, we show that enforcing these as hard constraints leads to an ill-conditioned regression problem due to the large number of constraints.
Instead, we implement soft constraints by modifying the cost function to be minimized.
This leads to the formulation of a multi-objective sparse regression problem where we simultaneously seek to minimize the error of the fit to the transients trajectories and to the data on attractors, while penalizing the number of terms in the model.
Our extension, demonstrated on several numerical examples, is more robust to noisy measurements and requires substantially less training data than the original SINDy method to correctly identify a parameterized dynamical system. 
\end{abstract}

\maketitle


\section{\label{sec:intro}Introduction}

The behavior of complex systems found across science and engineering often depends on quantities such as material or media properties, geometrical parameters, the magnitude of the physical scales involved, boundary conditions, and/or control inputs. The effect of parameter changes in nonlinear dynamics has been the focus of intense research for decades due to the emergence of bifurcations and pattern formation~\citep{cross1993rmp}. Moreover, reduced-order models that capture the parameterized dynamics of complex engineering systems~\citep{benner2015sr,rowley2017arfm} are a critical enabler for real-time applications~\citep{brunton2015amr,herrmann2020se,herrmann2022prf}, many-query optimization tasks~\citep{gunzburger1999ijnmf,giles2000ftc,herrmann2016ate,herrmann2018hmt}, and the development of digital twins~\citep{grieves2017springer,hartmann2018springer,niederer2021natcs}, which will have a profound impact on industries such as aerospace~\citep{kapteyn2021natcs,mcclellan2022ptrsa}, electromobility~\citep{ali2023as}, and biomedical~\cite{arzani2021jrsi}.

Traditionally, the governing equations of dynamical systems have been derived from first principles. However, although we are be able to build physics-based models for many complex systems, these models may be computationally expensive to simulate or may show poor predictive performance due to simplified physics assumptions. Furthermore, governing equations derived from first principles are not available in fields such as neuroscience, epidemiology, and finance. Fortunately, the unprecedented availability and quality of data has led to the emergence of data-driven model discovery — a paradigm in which we learn governing laws directly from measurements of the system of interest~\citep{databook}. In recent years, a plethora of approaches have been proposed in response to this trend, such as combined machine learning and data assimilation frameworks~\citep{gottwald2021physd,gottwald2021chaos}, deep learning approaches~\citep{pathak2018prl,vlachas2018prsa,chen2018neurips,raissi2019jcp,li2020arxiv,rackauckas2020arxiv,lu2021natmi}, divide-and-conquer strategies~\citep{udrescu2020sciadv}, symbolic regression via genetic algorithms~\citep{bongard2007pnas,schmidt2009science,cranmer2020neurips}, Gaussian processes~\citep{raissi2017jcp}, nonlinear autoregressive algorithms~\citep{narmaxbook}, the dynamic mode decomposition (DMD)~\citep{schmid2010jfm,rowley2009jfm,dmdbook} and its variants~\citep{schmid2022arfm,herrmann2021jfm,baddoo2022prsa,baddoo2023prsa}, Koopman theory~\citep{mezic2005nd,mezic2013arfm,brunton2022sr}, operator inference~\citep{peherstorfer2016cmame,qian2020physd}, regression onto spectral submanifolds~\citep{cenedese2022natcomms,cenedese2022ptrsa,axaas2023nd}, and the sparse identification of nonlinear dynamics (SINDy)~\citep{brunton2016pnas,rudy2017sciadv}. The latter method, which is the focus of the present work, has gained a considerable amount of traction because it follows a simple philosophy, is easy to implement, and is highly extensible~\citep{brunton2016pnas}.

SINDy builds parsimonious and interpretable models of dynamical systems from data by leveraging sparsity promoting regression techniques~\citep{tibshirani1996jrsssb} that balance model complexity and accuracy to avoid overfitting~\citep{brunton2016pnas}. Since its introduction, the method has been rapidly adopted to address challenges in a wide variety of fields, including active matter~\citep{joshi2022prl,golden2023sciadv}, air pollution~\citep{rubio2022jcpr}, astrophysics~\citep{pasquato2022taj},  biology~\citep{mangan2016ieeetmbmc,lagergren2020prsa}, chemistry~\citep{boninsegna2018jchp,hoffmann2019jchp,bhadriraju2020aiche,scheffold2021cche}, cognitive science~\citep{dale2018csr}, discrepancy modeling~\citep{kaheman2019arxiv,desilva2020fai,ebers2022arxiv}, epidemiology~\citep{jiang2021nd}, fluid dynamics~\citep{loiseau2018jfm_a,loiseau2018jfm_b,el2018aiaam,chang2019cg,deng2020jfm,loiseau2020dtcfd,fukami2021jfm,guan2021rsoc,deng2021jfm,callaham2022sciadv,callaham2022jfm,khoo2022jfm,sansica2022aiaaj,xiao2023pof}, manufacturing~\citep{subramanian2021mla}, heat transfer~\citep{ito2023jcnd}, numerical analysis~\citep{thaler2019jcp}, optics~\citep{sorokina2016oe}, plasma physics~\citep{dam2017pop,kaptanoglu2021pre,alves2022prr}, power systems~\citep{stankovic2020ieeepesgm,cai2022ieeetps}, structural dynamics~\citep{lai2019mssp,chatterjee2023mssp}, and turbulence closure modeling~\citep{beetham2020prf,schmelzer2020ftc,zanna2020grl,beetham2021jfm}. Furthermore, a large number of innovations that extend the capabilities of the original SINDy algorithm have been developed to discover partial differential equations (PDEs)~\citep{rudy2017sciadv,schaeffer2017prsa},  identify hybrid~\citep{mangan2019prsa} and stochastic~\citep{klimovskaia2016ploscb,bruckner2020prl,dai2020chaos,callaham2021prsa,huang2022sr} dynamical systems, incorporate control or exogenous inputs~\citep{brunton2016ifac,kaiser2018prsoc,fasel2021cdc}, express system dynamics with rational functions~\citep{mangan2016ieeetmbmc,kaheman2020prsa}, enforce physical constraints such as symmetries and conservation laws~\citep{loiseau2018jfm_a}, global stability~\citep{kaptanoglu2021prsa}, and dimensional consistency~\citep{bakarji2022natcs}, 
add uncertainty quantification and probabilistic forecasting capabilities~\citep{zhang2018prsa,ram2021arxiv,hirsh2022rsos,fasel2022prsa,gao2023arxiv}, allow for online model updating~\citep{quade2018chaos}, leverage smart sampling strategies~\citep{schaeffer2018siads,champion2019siads}, propose model selection criteria~\citep{mangan2017prsa,dong2023nd}, and dramatically improve the robustness of the method to noisy~\citep{schaeffer2017pre,gurevich2019chaos,messenger2021mms,messenger2021jcp,delahunt2022ieeea,kaheman2022mlst}, corrupt~\citep{tran2017mms,champion2020ieeea}, and incomplete~\citep{reinbold2020pre,reinbold2021natcomm,somacal2022pre,bakarji2023prsa} data. Many of these variants have been implemented in the open source software package PySINDy~\citep{desilva2020joss,kaptanoglu2022joss}, and some of them have been recently benchmarked~\citep{kaptanoglu2023nd} using the dysts standardized database of chaotic systems introduced by~\cite{gilpin2021neurips}. Another foundational contribution is that of~\citet{champion2019pnas}, who used autoencoder neural networks  along with SINDy to simultaneously discover latent coordinates and dynamics from high-dimensional data. Subsequently, several works have leveraged autoencoders in conjunction with SINDy for dimensionality reduction and identifying dynamics, respectively, in novel methods for data-driven low-order modeling for nonlinear PDEs~\citep{fukami2021jfm,fries2022cmame,callaham2022jfm,conti2023cmame}.

In the context of parameterized systems, the original SINDy article already introduced a basic method that is able to successfully identify parameterized dynamics, if provided with a large and noise-free data set~\citep{brunton2016pnas}. Recently,~\citet{meidani2023eswa} used integer programming to drastically improve the noise robustness of the method, albeit for the specific scenario where the model coefficients are known to be integers. The works of~\citet{schaeffer2017arxiv} and~\citet{rudy2019siads} used group sparsity to identify dynamical systems and PDEs, respectively, from data sampled at multiple parameter values. The group sparsity results in the identification of one model structure with model coefficients that are allowed to change as the parameter of the system varies, thus enabling model discovery in a scenario where the parameter is unknown, slowly drifting, or abruptly changing~\citep{schaeffer2017arxiv,rudy2019siads}. Consequently, it is important to point out that, this approach does not produce a symbolic expression for the model dependence on the parameter. In recent work, \citet{nicolau2023prr} developed a simple extension of the PDE-FIND algorithm~\citep{rudy2017sciadv} to identify parameterized PDEs by carefully building the library of function candidates to represent the parameter dependence. This extension was shown to successfully produce models capable of extrapolation to predict bifurcations~\citep{nicolau2023prr}. However, the authors observe that the main challenges to discover accurate models with this approach are that the method requires sufficiently informative trajectory data and that sufficiently many parameter values should be measured~\citep{nicolau2023prr}. Interestingly, they also remark that, according to their findings, trajectories with persistent dynamics provide more information than transient trajectories. Also recently,~\citet{conti2023cmame} developed a data-driven method to build nonlinear low-order models of parameterized PDEs. They leveraged the autoencoder + SINDy framework~\citep{champion2019pnas,fukami2021jfm,fries2022cmame,callaham2022jfm} to identify parameterized dynamics on latent coordinates in a low-dimensional embedding~\citep{conti2023cmame}, which allowed for computationally cheap numerical continuation of periodic solutions. Importantly, successful identification of parameterized dynamical systems with any of the methods described above requires a rich dataset containing measurements of multiple transient trajectories with samples distributed across parameter space. Reducing the amount of transient data required for parameterized model discovery is the main subject of this work.

\begin{figure}
    \centering
    \includegraphics[width=1\textwidth]{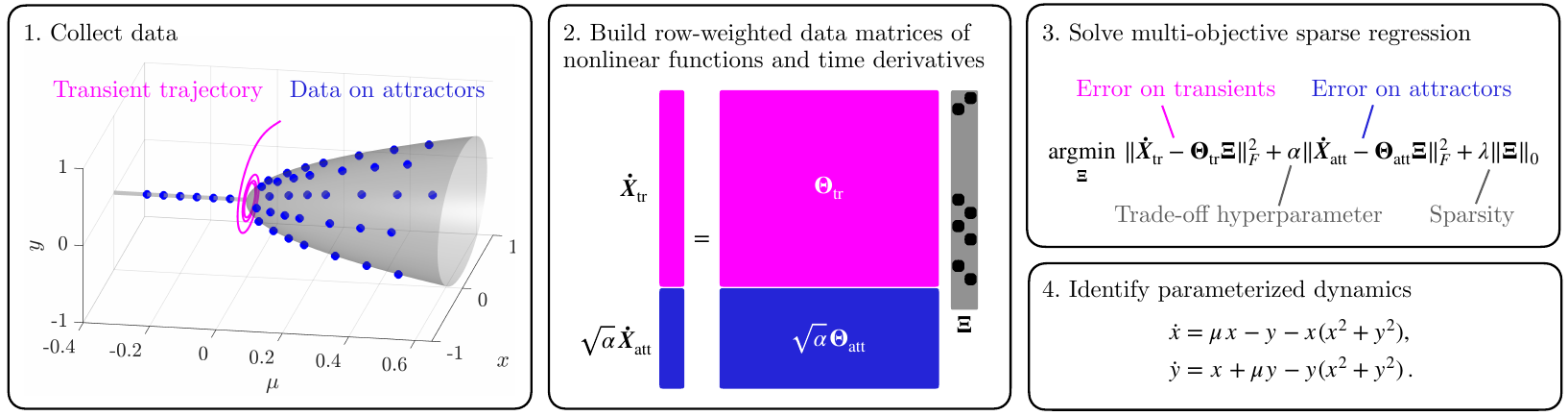}
    \caption{Schematic of the multi-objective SINDy framework for parameterized systems. The method balances measurements on and off attractors to accurately identify parameterized dynamics from sparse data. In the example, a Hopf normal form is correctly identified from a single transient trajectory complemented with a few measurements of stable fixed points and limit cycles distributed over parameter space.}
    \label{method}
\end{figure}

In an experimental setting, obtaining data on transient trajectories amounts to acquiring time-resolved measurements, whereas data on attractors is obtained from long-term response measurements. In addition, a stable fixed point can be accurately measured using non-time-resolved techniques, and a stable periodic orbit can be observed using stroboscopic measurements or phase averaging time-resolved data. In a numerical setting, we generate transient trajectories by time-marching through forward integration, whereas we may use iterative solvers to find fixed points and periodic orbits. Moreover, numerical continuation can be used to efficiently follow these solution branches through parameter space~\citep{dankowicz2013siam}. In short, different experimental techniques and numerical methods are typically used to obtain data on or off attractors. As a very general statement, measurements on attractors can typically be acquired at a lower computational cost (or with a simpler experimental setup) and often have a higher fidelity, compared to transient recordings of the same system. Furthermore, attractors and their bifurcations form the backbone of the state space of nonlinear dynamical systems and, because of this, we hypothesize that they provide more valuable data samples for parameterized system identification.

In this work, we develop a SINDy extension to identify parameterized dynamical systems using mostly data on attractors, and a minimum amount of data from transient trajectories. Our proposed approach, schematically shown in figure~\ref{method}, introduces a multi-objective sparse regression that balances the fit to both types of measurements to drastically improve the numerical conditioning of the regression.

The remainder of the paper is organized as follows. 
Theoretical background on SINDy for parameterized systems and SINDy with constraints is covered in~\S\ref{sec:background}. Our multi-objective SINDy framework is formulated in~\S\ref{sec:method}. Application of the method to several numerical examples is presented and discussed in~\S\ref{sec:results} and our conclusions are offered in \S\ref{sec:conclusions}.

\section{Background\label{sec:background}}

In this section we provide a brief review on the original SINDy formulation for parameterized systems presented in~\citet{brunton2016pnas} and the approach to incorporate constraints introduced in~\citet{loiseau2018jfm_a}.

\subsection{SINDy for parameterized systems}

The starting point is a parameterized nonlinear dynamical system with unknown dynamics that we want to identify. SINDy assumes that the dynamics can be expressed as a linear combination of a few nonlinear functions selected among a much larger, and judiciously designed, set of candidate functions, as follows
\begin{equation}
    \b{\dot{x}}=\b{f}(\b{x},\mu)\approx \b{\Xi}\T\b{\theta}(\b{x},\mu),
\end{equation}
where the overdot denotes time-differentiation, $\b{x}\in\mathbb{R}^n$ is the state of the system, $\mu\in \mathbb{R}$ is a parameter, $\b{f}: \mathbb{R}^n\times \mathbb{R} \rightarrow \mathbb{R}^n$ corresponds to the unknown dynamics, $\b{\theta} : \mathbb{R}^n\times \mathbb{R} \rightarrow \mathbb{R}^d$ contains the set of $d$ candidate scalar functions, and $\b{\Xi}\in\mathbb{R}^{d\times n}$ are the yet to be identified model coefficients.

Given a set of $m$ measurements of the state and parameter $\lbrace \b{x}_j,\mu_j\rbrace$ for $j= 1,\dots,m$, typically acquired from several transient trajectories generated for different parameter values, we may assemble the data matrices
\begin{equation*}
    \b{\dot{X}} = \left[\b{\dot{x}}_1 \ \b{\dot{x}}_2 \ \cdots \ \b{\dot{x}}_m\right]\T  \in\mathbb{R}^{m\times n},\quad \text{and } \quad
    \b{\Theta}=\left[\b{\theta}(\b{x}_1,\mu_1) \ \b{\theta}(\b{x}_2,\mu_2) \ \cdots \ \b{\theta}(\b{x}_m,\mu_m)\right]\T \in\mathbb{R}^{m\times d},
\end{equation*}
where  the time derivatives $\b{\dot{x}}_j=\dot{\b{x}}(t_j)$ may be approximated from sequential state data, via finite differences for example, if not directly measured. We may now formulate an optimization problem to identify $\b{\Xi}$, as follows
\begin{equation}
  \b{\Xi} = \underset{\b{\Xi'}}{\mathrm{argmin}} \ \| \b{\dot{X}}  -\b{\Theta} {\b{\Xi'}} \|_F^2 + \lambda \|\b{\Xi'}\|_0,\label{sindy_opt}
\end{equation}
where the first term in the cost function accounts for how well the model fits the data, and the second term, weighted by the hyperparameter $\lambda$, penalizes the amount of non-zero entries in the coefficient matrix, that is, the amount of terms in the identified model. Finding an exact solution to~\eqref{sindy_opt} requires performing a brute force search over all possible model structures, which is intractable when the number of candidate functions is large. In practice, sparse regression techniques, such as LASSO~\citep{tibshirani1996jrsssb} or the sequentially thresholded least squares (STLS) introduced in the original SINDy article~\citep{brunton2016pnas}, bypass the combinatorial search by only approximating the optimal solution to~\eqref{sindy_opt}.

Note that this formulation slightly differs from that presented for parameterized systems in the original article. In~\citep{brunton2016pnas} the state is simply augmented to contain the parameter as a new state variable governed by the trivial dynamics $\dot{\mu}=0$, which is included as an additional equation for the system. Although this may be convenient from a programming standpoint, because it allows reusing the same code developed for non-parameterized systems, it is not strictly necessary and the extra equation may be avoided altogether.

It is important to point out once again that robust identification of a parameterized dynamical system with SINDy typically requires a rich dataset comprised of multiple transient trajectories generated for multiple parameter values over the range of interest. In this work we develop a method to significantly reduce the amount of transient data required for successful identification of parameterized systems.

\subsection{SINDy with constraints}

The method introduced by~\citet{loiseau2018jfm_a} allows incorporating physical priors in the identification process by solving a constrained sparse regression. This is achieved by modifying the optimization problem~\eqref{sindy_opt} to
\begin{equation}
  \b{\xi} = \underset{\b{\xi'}}{\mathrm{argmin}} \ \| \b{\dot{\chi}}  -\b{\hat{\Theta}} {\b{\xi'}} \|_2^2 + \lambda \|\b{\xi'}\|_0, \quad \text{subject to} \quad \b{C}\b{\xi'} = \b{d},\label{csindy_opt}
\end{equation}
where $\b{\xi}$ and $\b{\dot{\chi}} \in \mathbb{R}^{nd}$ are vectorized forms of the $\b{\Xi}$ and $\b{\dot{X}}$ matrices, respectively, $\b{\hat{\Theta}}\in\mathbb{R}^{nm\times nd}$ is a block-diagonal matrix with $n$ copies of $\b{\Theta}$ as its blocks, and the matrix $\b{C}\in\mathbb{R}^{p\times nd}$ and vector $\b{d}\in\mathbb{R}^{p}$ enforce $p$ linear equality constraints. The optimization problem~\eqref{csindy_opt} is then approximately solved by iteratively solving a constrained least squares problem, finding small coefficients that satisfy $|\xi_i|<\lambda$, and augmenting $\b{C}$ and $\b{d}$ to incorporate the constraint $\b{e}_i\T \b{\xi} = \xi_i=0$ to promote sparsity.

The constrained least squares regression obtained for each iteration is a convex optimization problem, has an analytical solution that can be explicitly built for small-scale problems, and can be solved using readily available software for larger-scale ones~\citep{loiseau2018jfm_a}. Crucially, this solution exists if and only if the following conditions are satisfied:
\begin{enumerate}
    \item $\b{C}$ has linearly independent rows, and
    \item $\begin{bmatrix} \b{\hat{\Theta}}\\ \b{C} \end{bmatrix}$ has linearly independent columns.
\end{enumerate}
The first condition requires that $\b{C}$ is wide ($p\le nd$) meaning that there are fewer constraints than model coefficients. The second condition depends on both, $\b{\hat{\Theta}}$ and $\b{C}$, and it can be satisfied even if the function library has linearly dependent columns.

In addition to enforcing sparsity, constraints allow the incorporation of problem-specific knowledge of the system, such as fixed points, symmetries, and conservation laws~\citep{loiseau2018jfm_a}. The resulting models not only respect the known physics, but also generalize better and require less data for training. Moreover, in the context of identifying parameterized systems from as few measurements as possible across parameter space, SINDy with constraints seems like an attractive alternative to enforce parameter-dependent priors, such as fixed points and/or limit cycles. However, as we discuss in the next section, there are critical pitfalls to this approach.

\section{Proposed method\label{sec:method}}

In this section, we develop our SINDy extension to identify parameterized dynamical systems from limited data. We begin with our motivation to distinguish between data obtained on and off attractors, continue with the description of a first method that has limited use in practice, and then pivot to our proposed multi-objective sparse regression. 

\subsection{Data on and off attractors}

Acquisition of time-resolved measurements from transient trajectories and long-term response measurements from attractors often require different experimental setups. Even for the case of numerical data, we may solve for a stable steady or periodic state without having to perform a transient simulation. Furthermore, numerical continuation allows us to efficiently follow attractors in parameter space. This translates into a significant difference in cost, either economic or computational, between the acquisition of data on and off attractors. Moreover, attractors and their bifurcations provide a skeleton for the state space dynamics of parameterized systems. For these reasons, we hypothesize that an accurate representation of attractors is critical to successfully identify the correct model structure and, more so, it may even be sufficient to capture the parameter-dependence of the dynamics. This motivates the development of a method to identify parameterized dynamical systems relying as much as possible on data on attractors, and using as less data as possible from transient trajectories. Our approach to achieve this is to enforce the on-attractor dynamics as constraints.

We begin by distinguishing between data from transient trajectories $\lbrace \b{x}^\tr_j,\mu^\tr_j\rbrace$, for $j=1,\dots,m$, and data on attractors $\lbrace \b{x}^\att_j,\mu^\att_j\rbrace$, for $j=1,\dots,r$, leading to the definition of the following data matrices
\begin{align*}
    &\b{\dot{X}}_\tr = \left[\b{\dot{x}}^\tr_1 \ \b{\dot{x}}^\tr_2 \ \cdots \ \b{\dot{x}}^\tr_m\right]\T , \quad
    &\b{\Theta}_\tr=\left[\b{\theta}(\b{x}^\tr_1,\mu^\tr_1) \ \b{\theta}(\b{x}^\tr_2,\mu^\tr_2) \ \cdots \ \b{\theta}(\b{x}^\tr_m,\mu^\tr_m)\right]\T,\\
    &\b{\dot{X}}_\att = \left[\b{\dot{x}}^\att_1 \ \b{\dot{x}}^\att_2 \ \cdots \ \b{\dot{x}}^\att_r\right]\T ,\quad \text{and}
    &\b{\Theta}_\att =\left[\b{\theta}(\b{x}^\att_1,\mu^\att_1) \ \b{\theta}(\b{x}^\att_2,\mu^\att_2) \ \cdots \ \b{\theta}(\b{x}^\att_r,\mu^\att_r)\right]\T.
\end{align*}
Furthermore, if the attractor is a fixed point, then the time derivative is $\b{\dot{x}}^\att_j=\b{0}$, and, if it refers to a limit cycle, then it can be computed analytically from a Fourier series representation.

\subsection{The limitation of hard constraints}

As a first effort to formulate a method that leverages data on attractors across parameter space, we attempt using the constrained sparse identification developed in~\citep{loiseau2018jfm_a} and discussed in~\S\ref{sec:background}. We start from the optimization problem~\eqref{csindy_opt} using the data from the transient trajectories. In addition to the constraints $\b{C}\b{\xi} = \b{d}$, used to enforce the sparsity and other knowledge about the system, we include additional linear equality constraints $\b{\hat{\Theta}}_\att\b{\xi} = \b{\dot{\chi}}_\att$ to impose the on-attractor dynamics, as follows

\begin{equation}
  \b{\xi} = \underset{\b{\xi'}}{\mathrm{argmin}} \ \| \b{\dot{\chi}}_\tr  -\b{\hat{\Theta}}_\tr {\b{\xi'}} \|_2^2 + \lambda \|\b{\xi'}\|_0, \quad \text{subject to} \quad \b{\hat{\Theta}}_\att\b{\xi'} = \b{\dot{\chi}}_\att \ \text{ and } \ \b{C}\b{\xi'} = \b{d},
  \label{musindy_fail}
\end{equation}
where $\b{\hat{\Theta}}_\tr\in\mathbb{R}^{nm \times nd}$ and $\b{\hat{\Theta}}_\tr\in\mathbb{R}^{nr \times nd}$ are the block-diagonal matrices with $n$ copies of the corresponding data matrices for the function libraries, and $\b{\dot{\chi}}_\tr\in\mathbb{R}^{nm}$ and $\b{\dot{\chi}}_\att\in\mathbb{R}^{nr}$ are the vectorized data matrices with the time derivatives. Both sets of equality constraints can be combined to obtain
\begin{equation*}
    \b{\hat{\Theta}}_\att\b{\xi} = \b{\dot{\chi}}_\att \ \text{ and } \ \b{C}\b{\xi'} = \b{d} \qquad \Leftrightarrow \qquad 
    \b{\tilde{C}}\b{\xi}=\b{\tilde{d}}, \quad 
    \b{\tilde{C}} = \begin{bmatrix} \b{\hat{\Theta}}_\att\\ \b{C} \end{bmatrix}, \quad 
    \b{\tilde{d}} = \begin{bmatrix} \b{\dot{\chi}}_\att\\ \b{d} \end{bmatrix}.
\end{equation*}
Importantly, in order for a solution to the optimization problem to exist, we require that the matrix $\b{\tilde{C}}\in\mathbb{R}^{(p+nr)\times nd}$ is wide. Consequently, the dimensions of $\b{\tilde{C}}$ must satisfy $(p+nr) \le nd$, or equivalently,
\begin{equation}
    nr\le nd-p,\label{solv_cond}
\end{equation}
where $nr$ is the total number of constraints imposed to enforce the on-attractor dynamics, $nd$ is the total number of model coefficients being fit, and $p$ is the number of sparsity constraints, making $nd-p$ the amount of non-zero model coefficients. Therefore, the number of constraints that one can impose is, at most, the amount of active coefficients in the underlying dynamics. However, the assumption that the dynamics are governed by equations with only a few terms is at the core of the SINDy approach. Moreover, when attempting to enforce the on-attractor dynamics across a parameter range, even if it is a single fixed point measured at several parameter values, $nr$ grows with the number of samples $r$.

As a consequence, condition~\eqref{solv_cond} will rarely be satisfied in practice, making this approach extremely limited in its application to enforce the parameterized on-attractor dynamics. Therefore, we do not recommend using the method presented in this subsection for this purpose. Instead, we propose the approach presented next.

\subsection{Multi-objective sparse regression}

In view of the limited applicability of the method described above, we consider a different approach. Rather than 
incorporating the on-attractor dynamics as hard constraints, we include them in our cost function as soft constraints that merely penalize, as opposed to strictly enforcing, the model error on the attractor data. This leads to the formulation of a multi-objective sparse regression problem, as follows
\begin{equation}
  \b{\Xi} = \underset{\b{\Xi'}}{\mathrm{argmin}} \ 
  \underbrace{\| \b{\dot{X}}_\tr  -\b{\Theta}_\tr {\b{\Xi'}} \|_F^2}_{\textstyle{J_\tr(\b{\Xi'})}} 
  + \alpha \underbrace{\| \b{\dot{X}}_\att  -\b{\Theta}_\att {\b{\Xi'}} \|_F^2}_{\textstyle{J_\att(\b{\Xi'})}} 
  + \lambda \|\b{\Xi'}\|_0,\label{musindy}
\end{equation}
where the hyperparameter $\alpha$ weights how much we emphasize on minimizing the fit error on the attractors $J_\att(\b{\Xi})$ compared to the error on the transient trajectories $J_\tr(\b{\Xi})$. Note that, for a very small value of $\alpha$, the optimization will disregard the data on the attractors, whereas, for a very large $\alpha$, it will deemmphasize the transients. In fact, in the limit as $\alpha \rightarrow \infty$ we recover the constrained sparse regression~\eqref{musindy_fail}, which is ill posed if~\eqref{solv_cond} is not satisfied. Therefore, as we show in our numerical examples in \S\ref{sec:results}, there is typically a sweet spot for $\alpha$ that balances the fit to attractors and transients to find an accurate model across parameter space. As with any other hyperparameter, $\alpha$ may be selected via cross-validation or more sophisticated model selection criteria~\citep{mangan2017prsa,dong2023nd}.

Furthermore, the cost function in~\eqref{musindy} may be rewritten to yield
\begin{equation}
  \b{\Xi} = \underset{\b{\Xi'}}{\mathrm{argmin}} \
  \left\| \vphantom{\begin{bmatrix} . \\ . \end{bmatrix}} \right.
  \underbrace{\begin{bmatrix} \b{\dot{X}}_\tr\\ \sqrt{\alpha}\b{\dot{X}}_\att\end{bmatrix}}_{\textstyle{\b{\dot{X}}}} - 
  \underbrace{\begin{bmatrix} \b{\Theta}_\tr\\ \sqrt{\alpha}\b{\Theta}_\att\end{bmatrix}}_{\textstyle{\b{\Theta}}}
  {\b{\Xi'}} 
  \left. \vphantom{\begin{bmatrix} . \\ . \end{bmatrix}} \right\|_F^2
  + \lambda \|\b{\Xi'}\|_0,\label{musindy_stack}
\end{equation}

which has an equivalent structure to problem~\eqref{sindy_opt} and, for the case $\alpha=1$, we recover the original SINDy formulation~\citep{brunton2016pnas} using all the measurements stacked together. Thus, for a given $\alpha$, we can approximate the solution to~\eqref{musindy_stack} via STLS. In fact, we me reuse any available SINDy code after forming the row-weighted data matrices $\b{\dot{X}}$ and $\b{\Theta}$ in~\eqref{musindy_stack}. The pseudocode for our implementation is shown in algorithm~\ref{alg}, where we use Matlab's array slicing notation.

\begin{algorithm}[H]
\caption{Multi-objective SINDy for parameterized systems}\label{alg}
\begin{flushleft}
\textbf{Inputs:} \hspace{0.1 cm} data from transients $\b{\Theta}_\tr$ and $\b{\dot{X}}_\tr$\\
\hspace{1.35 cm} data on attractors $\b{\Theta}_\att$ and $\b{\dot{X}}_\att$\\
\hspace{1.35 cm} sparsification and trade-off hyperparameters $\lambda$ and $\alpha$\\
\textbf{Outputs:} SINDy model coefficients $\b{\Xi}$
\end{flushleft}
\begin{algorithmic}[1]
\State $\b{\dot{X}} \gets [\b{\dot{X}}_\tr; \sqrt{\alpha}\b{\dot{X}}_\att]$ \Comment{Weight and vertically concatenate the time derivatives.}
\State $\b{\Theta} \gets [\b{\Theta}_\tr; \sqrt{\alpha}\b{\Theta}_\att]$ \Comment{Weight and vertically concatenate the function libraries.}
\Function{STLS}{$\b{\dot{X}}$,$\b{\Theta}$,$\lambda$}
    \State $\b{v} \gets \max (\mathrm{abs}(\mathrm{columns}(\b{\Theta})))$ \Comment{Find the largest entry for each candidate function...}
    \State $\b{\Theta} \gets \b{\Theta}/\mathrm{diag}(\b{v})$ \Comment{...and normalize the columns of the function library.}
    \State $\b{\Xi} \gets \b{\Theta} \backslash \b{\dot{X}}$ \Comment{Initialize the coefficients with the least squares solution.}
    \For{$k=1,\dots,5$} \Comment{During each iteration...}
        \For{$i=1:n$} \Comment{...and for each state variable...}
            \State $\b{j}_{\mathrm{small}} \gets  \mathrm{abs}(\b{\Xi}(:,i)) < \lambda \max (\mathrm{abs}(\b{\Xi}(:,i)))$ \Comment{...find the entries with a small contribution to the dynamics...}
            \State $\b{\Xi}(\b{j}_{\mathrm{small}},i) \gets \b{0}$ \Comment{...and threshold.}
            \State $\b{j}_{\mathrm{big}} \gets \ \sim \b{j}_{\mathrm{small}}$ \Comment{Then find the non-zero entries...}
            \State $\b{\Xi}(\b{j}_{\mathrm{big}},i) \gets \b{\Theta}(:,\b{j}_{\mathrm{big}}) \backslash \b{\dot{X}}(:,i)$ \Comment{...and regress onto those terms.}
        \EndFor
    \EndFor
    \State $\b{\Xi} \gets \mathrm{diag}(\b{v}) \backslash \b{\Xi}$ \Comment{Rescale the rows of the coefficent matrix.}
\EndFunction
\end{algorithmic}
\end{algorithm}

There are two subtle details in the implementation of the STLS procedure that we use in algorithm~\ref{alg}.  First, because we have row-weighted data matrices, normalization of the columns of $\b{\Theta}$ becomes necessary to avoid having to scale the sparsity knob $\lambda$ with $\alpha$. Second, we introduce a variable-dependent thresholding criterion, where we retain the terms that most contribute to the dynamics of each state variable, rather than using one global thresholding criterion. We find that this becomes relevant when there is a significant variation in the size of the coefficients between the equations for different state variables.

A schematic of the application of the method to identify a Hopf normal form is shown in figure~\ref{method}. Our multi-objective SINDy framework successfully identifies the parameterized dynamics from data of a single transient trajectory along with stable fixed point and limit cycle data acquired for several parameter values. Due to its simplicity, the proposed approach may be trivially coupled with other powerful SINDy extensions, for example to identify PDEs~\citep{rudy2017sciadv}, add noise robustness~\citep{schaeffer2017prsa}, and quantify model uncertainty~\citep{fasel2022prsa}.

\section{Numerical examples and dataset\label{sec:data}}

\begin{figure}
    \centering
    \includegraphics[width=1\textwidth]{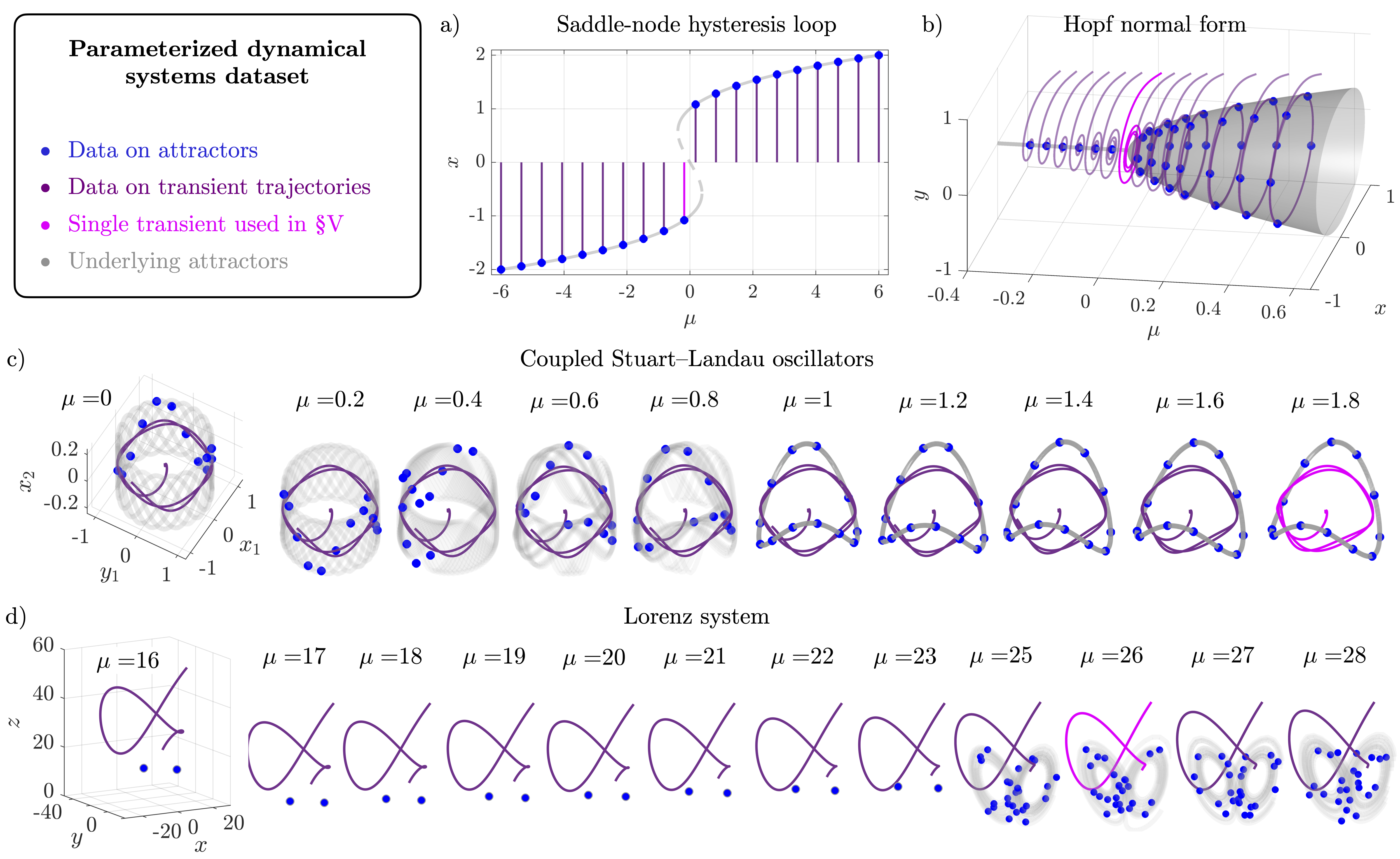}
    \caption{Numerical dataset of parameterized dynamical systems used to develop and compare methods for parameterized model discovery. For each system, off- and on-attractor measurements are acquired over a range of parameter values. The highlighted single transient trajectories are used for parameterized model discovery in~\S\ref{sec:results}. This dataset is available at github.com/ben-herrmann.}
    \label{fig:data}
\end{figure}

To demonstrate the application of multi-objective SINDy, we generate a dataset comprised of numerical measurements of four parameterized dynamical systems: a saddle-node hysteresis loop, the normal form of a supercritical Hopf bifurcation, a system of coupled Stuart–Landau oscillators, and the Lorenz equations. These systems are selected to feature all kinds of attractors; fixed points, limit cycles, quasiperiodic orbits, and strange attractors. For every system, we generate data by choosing an initial condition and numerically integrating \textit{long} trajectories for a range of parameter values. Numerical integration is carried out using an explicit $4-5$th order Runge–Kutta method implemented in the SciPy library in Python. Data points on and off attractors are registered, and a single transient trajectory, highlighted in figure~\ref{fig:data}, is flagged to be used for parameterized model discovery in~\S\ref{sec:results}. Data on transient trajectories are acquired at a fast sampling rate to appropriately resolve the time evolution. Time derivatives of the state variables are approximated with a $4$th order central finite difference scheme. Data on limit cycles, quasiperiodic and chaotic orbits are obtained using a burst sampling strategy~\citep{champion2019siads}, where short bursts of measurements are acquired at a fast sampling rate (the same used for transients) and separated by a much slower sampling rate to provide a better coverage of the attractor. Only the central data point within each burst is retained and their neighbors are used to approximate the corresponding time derivative and then discarded. Details on the selected dynamical systems and their data are presented below.

\subsection{Saddle-node hysteresis loop}

The dynamics for the saddle-node hysteresis loop are given by
\begin{equation}
    \dot{x}=\mu+x-x^3.
\end{equation}

Qualitatively similar models have been used to describe the dynamics of thermal runaway in batteries and insect outbreaks in forests~\citep{ludwig1978jae}. This system exhibits a pair of saddle-node bifurcations at $\mu=\pm 1/\sqrt{3}$, leading to two different branches of stable equilibria that coexist in the range $-1/\sqrt{3}<\mu<1/\sqrt{3}$.

Data for this system is acquired for parameter values starting at $\mu=-6$ and reaching $\mu=-0.182$ in steps of $0.65$, and also starting at $\mu=0.182$ and reaching $\mu=6$ in steps of $0.65$, for a total of $20$ sampled parameter values. For each $\mu$, transient trajectories are integrated starting from the initial condition $x_0 = 0$ and $1000$ data points are registered using a sampling time of $0.01$ time units. A single data point is registered on one of the stable equilibria branches for every parameter value. The negative branch is sampled for negative $\mu$ values, and the positive one for positive $\mu$ values. The single transient selected to test the proposed method corresponds to $\mu=-0.182$.

\subsection{Hopf normal form}

The governing equations for a supercritical Hopf normal form are
\begin{subequations}
\begin{align}
    \dot{x} &= \mu x -y -x(x^2+y^2), \\
    \dot{y} &= x +\mu y -y(x^2+y^2).
\end{align}
\end{subequations}

This is a well-known model for ubiquitous self-sustained oscillations appearing in fluid dynamics, chemistry and biology. The bifurcation occurs at $\mu=0$, where the fixed point $(0,0)\T$ looses its stability and trajectories end up in a stable limit cycle given by $\b{x}(t) = \sqrt{\mu}(\sin(t+\phi),\cos(t+\phi))\T$ with a phase $\phi$ that depends on the initial conditions.

Data for this system is acquired for parameter values starting at $\mu=-3$ and increasing up to $\mu=2$ in steps of $0.5$, and also at $\mu=0.3, \ 0.4,$ and $0.5$, for a total of $14$ sampled values. The value $\mu=0$ lands exactly at the bifurcation point so it is replaced by $\mu=0.01$. For each $\mu$, transient trajectories are integrated starting from the initial condition $\b{x}_0 = (2,0)\T$ and $2000$ data points are registered using a sampling time of $0.01$ time units. A single data point is registered on the stable equilibrium for each $\mu<0$, and, after the bifurcation, $5$ data points are acquired on the limit cycle sampled every $1.1$ time units for each $\mu>0$. The single transient selected to produce the results in~\S\ref{sec:results} corresponds to $\mu=0.01$.

\subsection{Coupled Stuart–Landau oscillators}

We consider two coupled Stuart–Landau oscillators governed by
\begin{subequations}
\begin{align}
    \dot{x}_1 &= -\omega_1 y_1 +\left(R_1^2 -(x_1^2+y_1^2)\right)x_1 +\mu(x_1 x_2 + y_1 y_2), \\
    \dot{y}_1 &= \omega_1 x_1 +\left(R_1^2 -(x_1^2+y_1^2)\right)y_1 +\mu(x_1 y_2 - y_1 x_2), \\
    \dot{x}_2 &= -\omega_2 y_2 +\left(R_2^2 -(x_2^2+y_2^2)\right)x_2, \\
    \dot{y}_2 &= \omega_2 x_2 +\left(R_2^2 -(x_2^2+y_2^2)\right)y_2.
\end{align}
\end{subequations}

A similar set of equations were recently proposed as a universal model for periodically forced turbulent oscillator flows~\citep{herrmann2020cp}. Here, the unforced frequencies and amplitudes of each oscillator are fixed to $R_1=1, \ R_2 = 0.2,\  \omega_1=1/\pi,$ and $\omega_2=1$, respectively. The only free parameter left is $\mu$ that represents the strength of the coupling between the oscillators. Note that the coupling is only in one direction, with the second oscillator driving the first one. When the second oscillator reaches its periodic orbit, this driving acts as periodic forcing on the first oscillator. Because the frequencies $\omega_1$ and $\omega_2$ are incommensurate, if the forcing is weak, the system evolves towards a quasiperiodic attractor. At $\mu=|\omega_2-2\omega_1|/(2 R_2)\approx 0.909$ the system undergoes a bifurcation and, for stronger forcing, a phenomenon known as $2:1$ synchronization emerges~\citep{herrmann2020cp}. What happens is that the first oscillator phase locks to half the frequency of the second one, leading to a periodic solution corresponding to a stable limit cycle that appears as the underlying attractor.

Data for this system is acquired for $10$ parameter values ranging from $\mu=0$ to $\mu=1.817$ in steps of $0.202$. For each $\mu$, transient trajectories are integrated starting from the initial condition $\b{x}_0 = (10^{-4},0,10^{-2},0)\T$ and $5020$ data points are registered using a sampling time of $0.01$ time units. After integrating for $2000$ time units, a trajectory is considered to be on the underlying attractor, and $12$ data points sampled every $1$ time units are registered for each $\mu$. The single transient selected to produce the results in~\S\ref{sec:results} corresponds to $\mu=1.817$.

\subsection{Lorenz system}

We consider the famous Lorenz equations
\begin{subequations}
\begin{align}
    \dot{x} &= \sigma (y-x),\\
    \dot{y} &= x(\mu-z)-y,\\
    \dot{z} &= xy-\beta z,
\end{align}
\end{subequations}
that represent a simple model of heat convection in the atmosphere~\citep{lorenz1963jas}.
We take $\sigma=10$, $\beta=8/3$, and $\mu$ takes the place of the parameter classically denoted as $\rho$. We focus in the range $16\le \mu \le 28$, where the system has $2$ stable equilibria at $(\pm \sqrt{\beta(\mu-1)}, \pm \sqrt{\beta(\mu-1)}, \mu-1)\T$ that become unstable after undergoing a subcritical Hopf bifurcation when increasing the parameter over $\mu\approx 24.74$, leading to the emergence of a chaotic or strange attractor.

Data for this system is acquired for $\mu=16$ up to $\mu=28$ in increments of $1$, for a total of $12$ sampled parameter values. For each $\mu$, transient trajectories are integrated starting from the initial condition $\b{x}_0 = (30,-30,50)\T$ and $500$ data points are registered using a sampling time of $0.001$ time units. For $\mu<24.74$, one data point is registered for only on of the stable fixed points. For values of the parameter after the bifurcation, we integrate \textit{long} trajectories and, after $70$ time units, we consider them to be on the chaotic attractor and register $25$ data points sampled every $0.2$ time units. The transient trajectory corresponding to $\mu=26$ is selected to test the multi-objective SINDy method in~\S\ref{sec:results}.

\section{Results and discussion\label{sec:results}}

\begin{figure}
    \centering
    \includegraphics[width=1\textwidth]{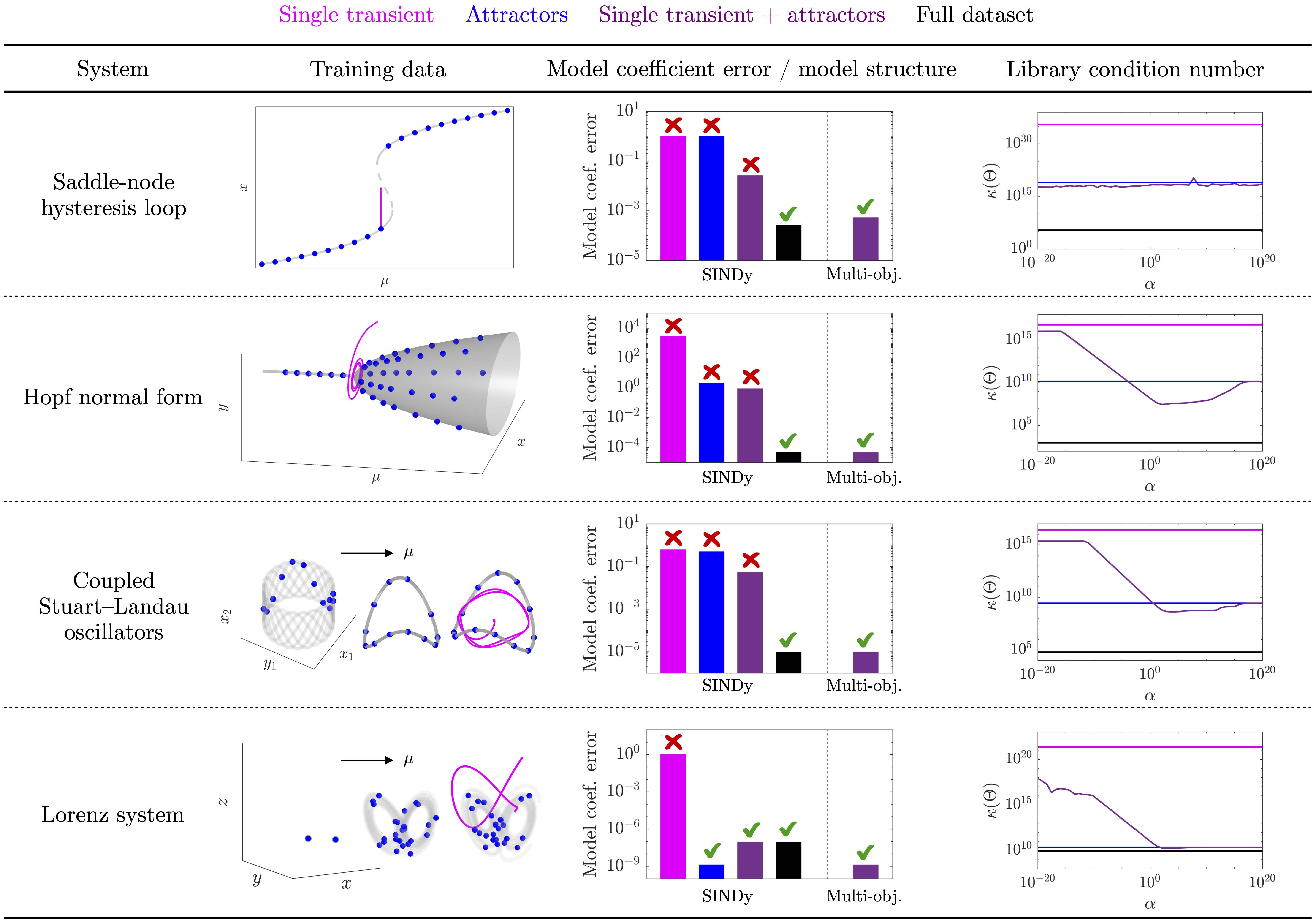}
    \caption{Multi-objective SINDy applied to a number of examples. Results are obtained using data from a single transient trajectory and from attractors measured over a parameter range. Performance at parameterized model discovery is assessed in terms of model coefficient error and the correctness of the identified model structure. Multi-objective SINDy models perform significantly better than standard SINDy models trained on the same data, and almost as well as standard SINDy models trained on a much larger dataset. This is explained because the row-weighting provided by the trade-off hyperparameter $\alpha$ drastically improves the numerical conditioning of the function library data matrices.}
    \label{fig:results}
\end{figure}

\begin{figure}
    \centering
    \includegraphics[width=\textwidth]{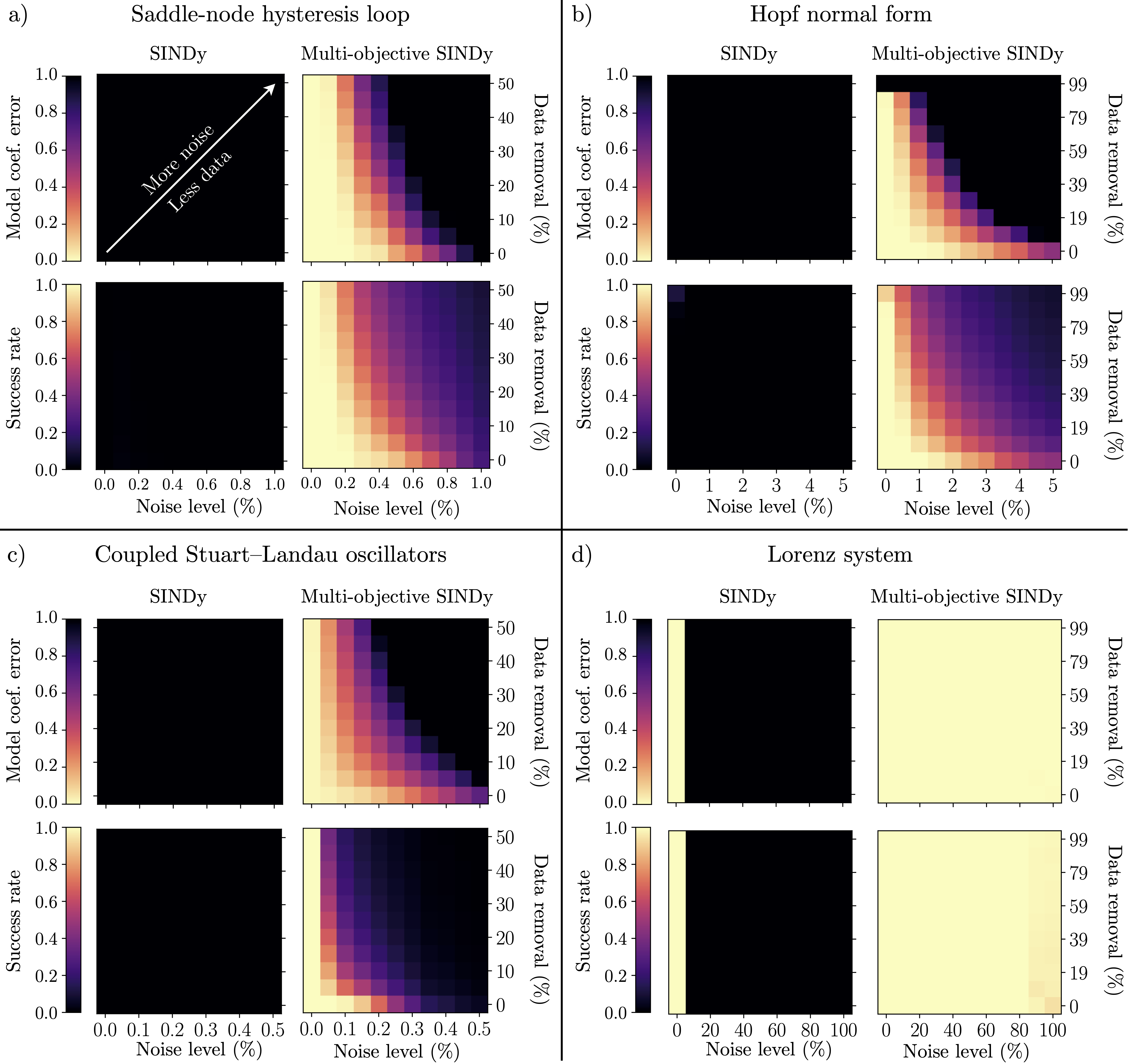}
    \caption{Robustness to noise and data length requirements for parameterized model discovery using data on attractors combined with sparse noisy data from a single transient trajectory over a range of parameters. Performance metrics are averaged over $500$ random realizations of measurement noise and data removal. (\textit{a}) Saddle-node hysteresis loop. (\textit{b}) Hopf normal form. (\textit{c}) Coupled Stuart–Landau oscillators. (\textit{d}) Lorenz system. Comparison between standard SINDy and multi-objective SINDy shows that the proposed method is able to identify more accurate models using less and more noisy data.}
    \label{fig:heatmaps}
\end{figure}

In this section, we demonstrate the application of multi-objective SINDy using the data for the systems described in the previous section. First, we apply our proposed method to identify the parameterized dynamics of these systems using a single transient trajectory and data on the corresponding attractors, shown in figure~\ref{fig:results}. We then investigate the noise robustness and data requirements of the method, shown in figure~\ref{fig:heatmaps}.

For each system, time derivatives are approximated via finite differences using a $4$th order central scheme. As mentioned in the previous section, time derivatives for data points on the attractors rely on the burst sampling strategy~\cite{champion2019siads}. For the function libraries, we use monomials of up to degrees $5, \ 4, \ 3$, and $4$ in the state variables and parameters for the saddle-node hysteresis loop, Hopf normal form, coupled Stuart-Landau oscillators, and Lorenz system, respectively.

\subsection{Parameterized system identification from a single transient trajectory}

The model identification performance of our method is compared to that of the standard SINDy algorithm from~\citep{brunton2016pnas} using a combination of data on the attractors and from a single transient trajectory, as shown in figure~\ref{fig:results}. Furthermore, to provide a more complete picture, we also present the performance of the standard SINDy algorithm using only the data on the attractors, only the data from the single transient trajectory, and the full dataset, shown in figure~\ref{fig:data}, containing multiple transient trajectories. Only noise-free numerical data is considered in this section. Models are identified for each system and dataset using various values for the hyperparameters $\lambda$ and $\alpha$. For the the sparsification knob $\lambda$ we do a sweep over a grid of $100$ log-spaced values ranging from $10^{-4}$ to $1$. When using multi-objective SINDy, the trade-off hyperparameter $\alpha$ is evaluated over a grid of $11$ log-spaced values between $10^{-20}$ and $10^{20}$. The hyperparameters that produce the best performing models for each case are shown in table~\ref{tab:lambdas}.
\begin{table}[h]
    \centering
    \begin{tabular}{c|c|c|c|c|c}
    \hline
       \multirow{2}{4em}{System}  & \multicolumn{4}{c|}{$\lambda$ (SINDy)} & $\lambda, \ \alpha$ (Multi-obj.) \\
         & Single transient & Attractors & Single trans.$+$attr. & Full dataset & Single trans.$+$attr. \\
         \hline
         Saddle-node hysteresis loop & $1.000$  &  $0.0001$  &  $0.0105$ &  $0.0095$  &  $0.0045, \ 10^{20}$\\
         Hopf normal form & $1.000$ &   $0.3944$  & $0.2719$  &  $0.0003$   & $0.0060, \ 10^4$ \\
         Coupled Stuart–Landau oscillators & $0.0001$ &  $0.0152$ &   $0.0034$ & $0.0001$ & $0.0087, \ 10^4$ \\
         Lorenz system & $0.0072$  & $0.0001$ &  $0.0001$ & $0.0001$ & $0.0001, \ 10^{20}$ \\
         \hline
    \end{tabular}
    \caption{Hyperparameters used to compare the performance of SINDy and multi-objective SINDy. For each system and training data, these values produce the lowest model coefficient error over a hyperparameter sweep.}
    \label{tab:lambdas}
\end{table}

The performance of the identified models is assessed in terms of model coefficient error defined as
\begin{equation}
    \text{Error} = \frac{\|\b{\Xi}-\b{\Xi}_{\mathrm{true}}\|_F}{\|\b{\Xi}_{\mathrm{true}}\|_F},
\end{equation}
where $\b{\Xi}$ are the identified model coefficients and $\b{\Xi}_{\mathrm{true}}$ are the true coefficients from the underlying system. In each case, we assess whether or not the correct model structure was identified based on the sparsity pattern of $\b{\Xi}$. Results shown in figure~\ref{fig:results} correspond to the best performing models in each case using the hyperparameters from table~\ref{tab:lambdas}.

Notice that, in general, SINDy is not able to successfully identify a parameterized dynamical system neither from a single transient trajectory nor from data on attractors only. This is expected because, in the first case we have not provided any information about the parameter dependence of the system. In the second case, identification fails because SINDy finds the model with the fewest number of terms to fit the data on the attractors, which is often a simpler model than the one that also explains the off-attractor dynamics. To provide some intuition, the simplest model that explains a fixed point measurement is given by the trivial dynamics $\b{\dot{x}}=\b{0}$, and the simplest model that explains observations from a limit cycle is given by a linear model with a center equilibrium.

The exception to this is the Lorenz system that has a chaotic attractor that provides enough information for successful for parameterized system identification. However, for a system where the chaotic attractor lies on a subspace of the full state space, measurements of transient trajectories would still be required to identify the off-attractor dynamics. Such an example of a chaotic system where data on the chaotic attractor is not enough for model discovery can be easily produced. For example, this would be the scenario for the Lorenz system if we added a fourth state variable $w$ with dynamics $\dot{w}=-w$. On the attractor we would only measure $w=0$ and therefore would not be able to identify the underlying dynamics.

Moreover, for the presented systems, standard SINDy is also not able to correctly identify the underlying dynamics when using the combined data from a single transient and attractors, and requires the full dataset containing multiple transient trajectories, as shown in figure~\ref{fig:results}. We find that just by weighting the data matrices, multi-objective SINDy correctly identifies the model structure using data from just a single transient trajectory and on the attractors. Furthermore, the performance in terms of model coefficient error is almost as good as that obtained by the standard SINDy algorithm using the full dataset.

Again, the Lorenz system is the exception to this general behavior. Interestingly, for this case, data on the attractors is sufficient for successful identification and, more so, including the transient trajectories has a detrimental effect on model performance. Nevertheless, with multi-objective SINDy we achieve a similar performance because the weighting by $\alpha$ diminishes the relative contribution of the transient data to the regression results.

To further understand why the proposed method performs so well in comparison the the standard approach, we compute the condition numbers of the function library data matrices $\kappa(\b{\Theta})$ obtained using different parts of the dataset. In the last column of figure~\ref{fig:results}, we compare the condition numbers of the data matrices formed using the full dataset $\b{\Theta}_{\mathrm{full}}$ (in black), data from a single transient trajectory $\b{\Theta}_{\tr}$ (in magenta), only on the attractors $\b{\Theta}_{\att}$ (in blue), and the combination of the latter two weighted by the trade-off hyperparameter $\b{\Theta}_{\alpha}=\left[\b{\Theta}_{\tr}; \sqrt{\alpha}\b{\Theta}_{\att}\right]$ (in purple) as a function of $\alpha$. The last is the data matrix involved in the multi-objective sparse regression problem from equation~\eqref{musindy_stack} solved by the proposed method, and corresponds to the data matrix used in standard SINDy when $\alpha=1$. The plots for the Hopf normal form and coupled Stuart–Landau oscillators in figure~\ref{fig:results} clearly show that in the limit as $\alpha \rightarrow -\infty$ the condition number $\kappa(\b{\Theta}_{\alpha})\rightarrow \kappa(\b{\Theta}_{\tr})$, in the limit as $\alpha \rightarrow \infty$ then $\kappa(\b{\Theta}_{\alpha})\rightarrow \kappa(\b{\Theta}_{\att})$, and that there is a minimum in between that gets much closer to $\kappa(\b{\Theta}_{\mathrm{full}})$ than either extrema. Importantly, this minimum in the condition number does not necessarily occur at $\alpha=1$, and its position will depend on the relative amount of data points from each type (on- or off-attractor dynamics) and the data values themselves. As a consequence, the multi-objective approach significantly improves the conditioning of the regression problem by introducing a value of $\alpha\ne 1$ to better balance the information provided by the different types of measurements. For the Lorenz system we find that the library evaluated on data on the attractors has a condition number that is almost as low as that of the library evaluated in the full dataset and is significantly better conditioned than the library evaluated on a single transient. Moreover, for values of $\alpha<1$, increasing the relative contribution of the transient trajectories, we find that the conditioning drastically worsens, explaining the detrimental effect of including the transient measurements on model performance. Results shown here use the values of $\alpha$, shown in table~\ref{tab:lambdas}, that produce the best performing model in each case with the aim of highlighting the potential of multi-objective SINDy. However, in a real system identification scenario, we propose selecting the value of $\alpha$ via cross-validation or, as with any other hyperparameter, it can also be determined using other model selection criteria~\citep{mangan2017prsa,dong2023nd}.

\subsection{Robustness to noise and data requirements}

Improving the conditioning of the regression should provide robustness to the identification procedure when dealing with sparse and noisy measurements. Following the methodology of~\citet{fasel2022prsa}, we demonstrate this robustness by adding random measurement noise to our data and assessing the average performance over multiple realizations for a range of noise intensities and data lengths. Specifically, for each of the systems presented in~\S\ref{sec:data}, we use the data from a single transient trajectory and on the attractors for various values of the parameter. We add Gaussian white noise with zero mean and variance $\sigma^2$ to the transient only for a range of noise levels, defined as $\sigma$ normalized by the root mean square of the transient trajectory data. In addition, we randomly discard measurements of this same transient trajectory and retain a certain percentage of the original number of data points. The model coefficient error and success rate ($1$ if the correct model structure is identified, $0$ if not) are averaged over $500$ random realizations of noise and of discarded data points for each noise level and percentage of data removal, as shown in figure~\ref{fig:heatmaps}. For each system, the hyperparameters are fixed for all noise levels and data removal percentages to the values shown in table~\ref{tab:lambdas_noise}. These values are lightly tuned to make the trends in figure~\ref{fig:heatmaps} clear enough to support our main conclusions. However, we remark that these values have not been optimized in any way and that the performance of multi-objective SINDy could probably be improved even more. Furthermore, in a realistic scenario, the hyperparameters would be selected by cross-validating against data with the same noise level and sparsity as the training data, leading to hyperparameters that are specifically tuned for the characteristics of the available measurements.

\begin{table}[h]
    \centering
    \begin{tabular}{c|c|c}
    \hline
        System  & $\lambda$ (SINDy) & $\lambda, \ \alpha$ (Multi-objective SINDy) \\
         \hline
         Saddle-node hysteresis loop &  $0.15$  &  $0.15, \ 10^{14}$\\
         Hopf normal form &  $0.06$   & $0.06, \ 10^14$ \\
         Coupled Stuart–Landau oscillators & $0.0087$ & $0.0087, \ 10^10$ \\
         Lorenz system & $0.01$ & $0.01, \ 10^{14}$ \\
         \hline
    \end{tabular}
    \caption{Hyperparameters used to compare the robustness to noise and data requirements of SINDy and multi-objective SINDy. Models are trained using data on attractors combined with sparse noisy data from a single transient trajectory}
    \label{tab:lambdas_noise}
\end{table}

As mentioned in the previous subsection, using a single transient trajectory and the attractors, SINDy is not able to identify the parameterized dynamics for the saddle-node hysteresis loop, Hopf normal form or coupled Stuart–Landau oscillators, even when the data is clean and full length, that is with a noise level of $0\%$ and $0\%$ data removal. With multi-objective SINDy, not only can we identify these systems from clean full-length data, but the method significantly moves the front of successful identification towards the top-right corner of high intensity noise and sparse measurements, as shown in figure~\ref{fig:heatmaps}. Furthermore, although SINDy is able to identify the parameterized Lorenz dynamics without noise, multi-objective SINDy is more robust to noise. As discussed in the previous section for the Lorenz example, including the data on the transient trajectory, as opposed to just using data on the attractors, is actually detrimental for system identification performance.
Notably, this effect is not seen for multi-objective SINDy, that deemphasizes the data from the transient trajectory via the row-weighting of the data matrices to achieve better performance.

\section{Conclusions\label{sec:conclusions}}

Measurements of the state of a dynamical system evolving on an attractor typically have a lower acquisition cost and are more informative for parameterized system identification than measurements from transient trajectories. In this work, we developed a method to identify parameterized dynamical systems relying as much as possible on data on attractors, and using as less data as possible from transient trajectories. The approach is an extension of the SINDy method~\citep{brunton2016pnas} that introduces a multi-objective cost function to balance the fit to data on and off attractors, leading to a significant improvement in the numerical conditioning of the underlying regression problem. We show that parameterized nonlinear dynamics can be learned from a single transient trajectory and data on attractors acquired over a range of parameter values. The new multi-objective SINDy is more robust to noisy measurements and requires less data for successful model discovery than the standard SINDy algorithm.

Moreover, the developed method is very simple to implement and can be easily coupled with other available SINDy extensions. Furthermore, the key idea behind multi-objective SINDy can be adopted for problems beyond parameterized dynamical systems, in other contexts where data from multiple sources can be combined for model discovery. Future work considers using the approach in conjunction with dimensionality reduction techniques to develop parameterized low-order models for parameterized PDEs. An interesting challenge that we foresee following that path is that, not only the dynamics may depend on the parameter, but also the low dimensional embedding.

\begin{acknowledgments}
We gratefully acknowledge S. Brunton, U. Fasel, M. Bahamondes, D. Delgado, and N. Torres for helpful comments and insightful discussions.
This work was funded by U. of Chile internal grant U-Inicia-003/21 and ANID Fondecyt 11220465.
\end{acknowledgments}



\providecommand{\noopsort}[1]{}\providecommand{\singleletter}[1]{#1}%

\end{document}